\documentclass[11pt]{article}

\usepackage{amsmath,amssymb,amsfonts,enumitem,amsthm,accents,pstricks,pstricks-add,mathrsfs}

\newtheorem{thm}{Theorem}
\newtheorem{cnj}[thm]{Conjecture}
\newtheorem{ques}[thm]{Question}

\theoremstyle{definition}
\newtheorem{dfn}[thm]{Definition}
\newtheorem{rmk}[thm]{Remark}
\newtheorem{eg}[thm]{Example}

\def\BE#1{\begin{equation}\label{#1}}
\def\EE{\end{equation}}
\def\eref#1{(\ref{#1})}

\def\bIt{\begin{itemize}[leftmargin=*]}
\def\eIt{\end{itemize}}

\def\tn#1{\textnormal{#1}}
\def\sf#1{\textsf{#1}}

\def\wh#1{\widehat{#1}}
\def\wt#1{\widetilde{#1}}

\def\lra{\longrightarrow}
\def\lhra{\ensuremath{\lhook\joinrel\relbar\joinrel\rightarrow}}

\def\C{\mathbb C}
\def\cC{\mathcal C}
\def\sC{\mathscr C}
\def\nd{\tn{d}}
\def\fD{\mathfrak D}
\def\ne{\tn{e}}
\def\fI{\mathfrak i}
\def\fj{\mathfrak j}
\def\cJ{\mathcal J}
\def\cN{\mathcal N}
\def\cO{\mathcal O}
\def\P{\mathbb P}
\def\R{\mathbb R}
\def\cR{\mathcal R}
\def\Q{\mathbb Q}
\def\cS{\mathcal S}
\def\ft{\mathfrak t}
\def\Z{\mathbb Z}
\def\cZ{\mathcal Z}

\def\al{\alpha}
\def\la{\lambda}
\def\na{\nabla}
\def\om{\omega}
\def\th{\theta}
\def\vp{\varpi}
\def\ze{\zeta}

\def\De{\Delta}
\def\Ga{\Gamma}

\def\Om{\Omega}
\def\Si{\Sigma}

\def\AK{\tn{AK}}
\def\Aux{\tn{Aux}}
\def\Dom{\tn{Dom}}

\def\id{\tn{id}}
\def\Im{\tn{Im}}
\def\sing{\tn{sing}}
\def\Symp{\tn{Symp}}

\def\bu{\bullet}
\def\dbar{\bar\partial}
\def\eset{\emptyset}
\def\hb{\hbar}
\def\prt{\partial}

\begin{document}
\title{Singularities and Semistable Degenerations 
for Symplectic Topology}
\author{Mohammad F.~Tehrani, Mark McLean, and 
Aleksey Zinger\thanks{Partially supported by NSF grant 1500875}}
\date{\today}
\maketitle

\begin{abstract}
\noindent
We overview our work \cite{FMZDiv,FMZSum,FMZNC,FZCut,FZRelation} defining and studying 
normal crossings varieties and subvarieties in symplectic topology. 
This work answers a question of Gromov on the feasibility of introducing 
singular (sub)varieties into symplectic topology in the case of normal crossings singularities.
It also provides a necessary and sufficient condition for smoothing 
normal crossings symplectic varieties.
In addition, we explain some connections with other areas of mathematics 
and discuss a few directions for further research.
\end{abstract}

\section{Introduction}
\label{intro_sec}

\noindent
Singularities are ubiquitous in algebraic geometry. 
Even if one wishes to study only smooth varieties, 
one is often forced to investigate certain singular varieties as well.
This is especially true in enumerative geometry, moduli theory, and 
the minimal model program.
Degeneration techniques are particular instances where such singularities arise.
Suppose, for example, that one wishes to study a particular deformation-invariant property 
of a smooth variety~$X$ such as Gromov-Witten invariants.
One might then try to degenerate~$X$ through a flat one-parameter family to 
a simpler, but possibly singular, variety~$X_0$ and study~$X_0$ instead.
For instance, when one wishes to compute Gromov-Witten invariants,
the degeneration formulas of \cite{Jun,BP3} can be particularly useful.
In mirror symmetry, such ideas have been suggested in \cite{Morrison} and are
central to the Gross-Siebert program~\cite{GS0}.\\

\noindent
Gromov \cite[p343]{Gr} asked if there was an appropriate
notion of singular variety in symplectic topology and when such varieties could be smoothed.
Introducing singularities and degeneration techniques into symplectic topology 
would be extremely useful for the following reasons.
First, symplectic manifolds are significantly more flexible than algebraic varieties,
since one can apply Moser style arguments and sometimes even an h-principle~\cite{Gr}.
Smoothing or degenerating such symplectic varieties should then become 
a problem of a topological nature, and 
therefore one would not need to rely on subtle analytic invariants.
Second, these techniques could be used to study a much larger class of symplectic manifolds
(not just those coming from K\"{a}hler geometry).
For instance, they could be useful in mirror symmetry for non-K\"{a}hler
symplectic manifolds; see for example \cite[Section 2.3]{TsengYau}.
They could also be used to construct interesting examples of symplectic manifolds
by smoothing singular ones, as has been done to great effect in~\cite{Gf}.\\

\noindent
This note is an overview of our work \cite{FMZDiv,FMZSum,FMZNC,FZCut,FZRelation},
which introduces and studies symplectic topology notions of 
normal crossings (or {\sf NC}) divisor and variety.
In Section~\ref{SCdfn_sec}, we explain what these notions~are.
In Section~\ref{RegDfn_sec}, 
we describe geometric notions of \sf{regularization} for NC symplectic divisors and varieties,
which is basically a ``nice" neighborhood of the singular locus.
Every NC symplectic divisor/variety is deformation equivalent to one with a regularization.
In fact, we propose~to\\

\begin{minipage}{6in}
{\it view a symplectic sub(variety) as a deformation equivalence class of objects,
not as an individual object}.\\
\end{minipage}

\noindent
In Section~\ref{SympSum_sec},
we explain our result on the smoothability of NC symplectic varieties.
It provides a purely topological necessary {\it and} sufficient condition 
for such a variety to be smoothable.
In Section~\ref{SympCut_sec},
we explain how to use certain local Hamiltonian torus actions to degenerate 
a smooth symplectic manifold into an NC symplectic variety.
Finally, in Section \ref{FurRes_sec}, 
we discuss connections with other areas of mathematics and directions for further developments.\\

\noindent
We would like to thank A.~Cannas da Silva, K.~Fukaya, E.~Ionel, E.~Lerman, B.~Parker, 
A.~Pires, H.~Ruddat, D.~Stapleton, J.~Starr, and D.~Sullivan for enlightening discussions
on various aspects of our work \cite{FMZDiv,FMZSum,FMZNC,FZCut,FZRelation}.

\section{Normal crossings (sub)varieties} 
\label{SCdfn_sec}

\noindent
A \sf{symplectic manifold} is a manifold $X$ of an even real dimension~$2n$ 
together with a closed nondegenerate 2-form~$\om$, i.e.
$\nd\om\!=\!0$ and $\om^n|_x\!\neq\!0$ for every $x\!\in\!X$.
In particular, $\om^n$ is a volume/orientation form on~$X$. 
An \sf{almost complex structure} $J$ on~$X$ is a vector bundle endomorphism 
$J$ of~$TX$ covering~$\id_X$ such that $J^2\!=\!-\id_{TX}$.  
Every complex (holomorphic) structure on~$X$ determines an almost complex structure on~$X$, 
but the converse is not true if $n\!>\!1$. 
The \sf{Nijenhuis $(2,1)$-tensor},
$$N_J(\xi,\ze)=\frac14\big([\xi,\ze]+J[\xi,J\ze]+J[J\xi,\ze]-[J\xi,J\ze]\big) 
\in\Ga(X;TX) \qquad \forall \xi,\ze\in\Ga(X;TX),$$ 
is the obstruction to the \sf{integrability} of~$J$, 
i.e.~$J$ arises from a holomorphic structure on~$X$ if and only if $N_J\!\equiv\! 0$; 
see \cite{NN}.\\

\noindent
Every symplectic manifold $(X,\om)$ can be equipped with an \sf{$\om$-compatible} 
almost complex structure in the sense that $\om(\cdot,J\cdot)$ is a metric. 
The space $\cJ_{\om}(X)$ of $\om$-compatible almost complex structures 
is infinite-dimensional and contractible. 
Such a triple $(X,\om,J)$ is called an \sf{almost K\"ahler} manifold.  
A \sf{K\"ahler} manifold is an almost K\"ahler manifold such that 
the almost complex structure is integrable. 
The category of K\"ahler manifolds includes smooth complex projective varieties, 
i.e.~smooth varieties cut out by polynomial equations in a complex projective space~$\P^N$.\\

\noindent
In contrast to algebraic/K\"ahler geometry, symplectic topology
is significantly more ``flexible".
For instance, symplectic manifolds have no local invariants and smooth families 
of symplectic manifolds with cohomologous symplectic forms are symplectomorphic.
We are interested in finding out which algebraic structures can be generalized 
to useful structures in the symplectic topology category.\\

\noindent
A \sf{symplectic submanifold} of a symplectic manifold $(X,\om)$ 
is a submanifold~$V$ of~$X$ such that $\om|_V$ is a symplectic form.
Symplectic submanifolds are the analogues of smooth subvarieties in (complex) algebraic geometry. 
For example, a \sf{smooth symplectic divisor} is a symplectic submanifold of real codimension~2 
(or complex codimension~1) and a 
\sf{smooth symplectic curve} is a symplectic submanifold of real dimension~2. 
The \sf{normal bundle} 
\BE{cNXVsymp_e}
\pi_{\cN_XV}\!:
\cN_XV\equiv \frac{TX|_V}{TV}\approx TV^{\om}
\equiv \big\{v\!\in\!T_xX\!:\,x\!\in\!V,\,\om(v,w)\!=\!0~\forall\,w\!\in\!T_xV\big\}
\lra V\EE
of a symplectic submanifold $V$ of $(X,\om)$ inherits a fiberwise symplectic 
form~$\om|_{\cN_XV}$ from~$\om$. 
An \sf{$\om|_{\cN_XV}$-compatible Hermitian structure} on $\cN_XV$ is a triple $(\fI,\rho,\na)$,
where $\fI$ is an $\om|_{\cN_XV}$-compatible complex structure on~$\cN_XV$, 
$\rho$ is a Hermitian metric with 
$$\rho_{\R}(\cdot,\cdot)=\om|_{\cN_XV}(\cdot,\fI\cdot),$$ 
and $\na$ is a Hermitian connection compatible with $(\fI,\rho)$.
The space of $\om|_{\cN_XV}$-compatible Hermitian structures is non-empty and contractible.\\

\noindent
An $\om|_{\cN_XV}$-compatible Hermitian structure $(\fI,\rho,\na)$ as above
determines a 1-form $\al_{\na}$ on $\cN_XV\!-\!V$ 
whose restriction to each fiber $\cN_XV|_x\!-\{x\}\!\cong\!\C^*$ 
is the 1-form $\nd\th$ with respect to the polar coordinates $(r,\th)$ on~$\C$. 
The closed 2-form 
\BE{OmtoHatOm_e}\wh\om\equiv\pi^*(\om|_V)+\frac{1}{2}\nd(\rho\al_{\na}) \in \Om^2(\cN_XV) \EE
is well-defined, is nondegenerate in a small neighborhood of~$V$, and
restricts to the standard symplectic form $\nd(r^2\nd\th)$ on each fiber.
By the Symplectic Neighborhood Theorem \cite[Theorem~3.30]{MS1},  
there exists an identification (called an \sf{$\om$-regularization}, in what follows)
\BE{Reg1_e}
\Psi\!:\cN'_{X}V\lra X, \qquad \Psi|_{V}=\id_{V},\qquad 
\nd\Psi|_{V}=\id,\EE
of a small neighborhood $\cN'_{X}V$ of $V$ in $\cN_XV$
with a neighborhood of $V$ in $X$ such that $\Psi^*\om\!=\!\wh\om$.
Regularizations are useful for applications, such as the symplectic sum construction
of \cite{Gf,MW}.
They also ensure the existence of almost complex structures~$J$ on~$X$ that are ``nice"
along~$V$.
These are in turn useful for constructing relative Gromov-Witten invariants of~$(X,V)$,
for example.\\

\noindent
In the 1980s, Gromov combined the rigidity of algebraic geometry with the flexibility of the smooth category and initiated the use of \sf{$J$-holomorphic maps}
from Riemann surfaces $(\Si,\fj)$ into $(X,J)$,
\BE{JHoloMaps_e}
u\!: (\Si,\fj)\lra(X,J),\qquad 
\dbar u\!\equiv\! \frac{1}{2}\big(\nd u + J\nd u\circ\fj\big)=0,\EE
as a generalization of holomorphic maps.
The singularities of the image of a $J$-holomorphic map~$u$ are locally 
the same as the algebraic ones; see \cite[Appendix~E]{MS2}.
If $J$ is $\om$-compatible, then
the smooth locus of~$u$ is a symplectic submanifold of~$(X,\om)$. 
A singular symplectic variety in complex dimension~1 can thus be defined 
as a subset of $X$ that can be realized as the image of a $J$-holomorphic map. 
The spaces of $J$-holomorphic maps have a nice deformation theory,
which makes it possible to study these objects in families;
see \cite[Chapter~3]{MS2} and \cite[Section~3]{LT}, for example.
The idea of studying $J$-holomorphic maps from higher-dimensional domains 
is not as promising because the Cauchy-Riemann equation~\eref{JHoloMaps_e}
is over-determined if the dimension of~$\Si$ is greater than~2.\\

\noindent
In parallel with his introduction of $J$-holomorphic curve techniques into symplectic topology, Gromov asked about the feasibility of introducing notions of singular (sub-)varieties of higher dimension suitable for this field; see \cite[343]{Gr}. By the last paragraph, the idea of defining such singular objects as images of $J$-holomorphic maps is not promising and we should consider an intrinsic approach. Nevertheless, we still require that for such a singular object $V$ in $X$, the space (or a nice subspace) of almost complex structures on $X$ compatible 
with $V$ to be non-empty and ``manageable".\\

\noindent 
In algebraic geometry, \sf{divisors}, i.e.~subvarieties of codimension~1, 
are dual objects to curves and have long been of particular importance. 
On the symplectic side, smooth symplectic divisors appear in different contexts
such as in relation with complex line bundles~\cite{Donaldson},
symplectic sum constructions \cite{Gf,MW}, 
relative Gromov-Witten theory and degeneration 
formulas for Gromov-Witten invariants \cite{Tian,LR,IPrel,BP}, 
affine symplectic geometry \cite{MAff,McLean}, and homological mirror symmetry~\cite{Sheridan}. 
The following question is thus one of the most important specializations of Gromov's inquiry:\\

\begin{minipage}{6in}
{\it Can one define a soft notion of (singular) symplectic divisor 
that only involves soft intrinsic symplectic data, 
but at the same time is compatible with rigid auxiliary almost K\"ahler 
data needed for making such a notion useful?}\\
\end{minipage}

\vspace{.1in}

\noindent
NC divisors/varieties are the most basic and important classes of singular objects
in complex algebraic (or K\"ahler) geometry.
An \sf{NC divisor} in a smooth variety~$X$ is a subvariety~$V$ 
locally defined by an equation of the form
\BE{Local1_e} z_1\cdots z_k = 0\EE
in a holomorphic coordinate chart $(z_1,\ldots,z_n)$ on~$X$.
A \sf{simple normal crossings} (or {\sf{SC}) \sf{divisor} 
is a global transverse union of smooth divisors, i.e. 
$$V=\bigcup_{i\in S}\!V_i \subset X.$$
An~\sf{NC variety} of complex dimension $n$ is a variety $X_{\eset}$ 
that can be locally embedded as an NC divisor in~$\C^{n+1}$. 
In other words, every sufficiently small open set~$U$ in  $X_{\eset}$ 
can be written~as 
$$U=\bigg(\bigsqcup_{i\in S}\!U_i\bigg)\!\Big/\!\!\!\sim,
\qquad  U_{ij}\approx U_{ji}\quad \forall~i,j\!\in\!S,~i\!\neq\!j,$$
where $\{U_{ij}\}_{j\in S-i}$ is an SC divisor in a smooth component~$U_i$ of~$U$. 
A \sf{simple normal crossings} (or~\sf{SC}) \sf{variety} 
is a global transverse union of smooth varieties $\{X_i\}_{i\in S}$ 
along SC divisors $\{X_{ij}\}_{j\in S-i}$ in~$X_i$,~i.e.
$$X_{\eset}=\bigg(\bigsqcup_{i\in S}\!X_i\bigg)\!\Big/\!\!\!\sim,\qquad  
X_{ij}\approx X_{ji}\quad \forall~i,j\!\in\!S,~i\!\neq\!j.$$
A 3-fold SC variety is shown in Figure~\ref{3conf_fig}.\\

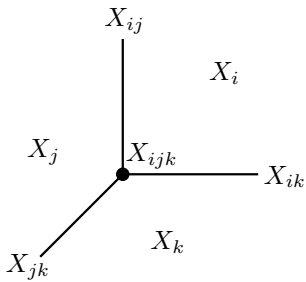
\begin{figure}
\begin{pspicture}(-3,-2)(11,.5)
\psset{unit=.3cm}
\psline[linewidth=.1](15,-2)(21,-2)\psline[linewidth=.1](15,-2)(15,4)
\psline[linewidth=.1](15,-2)(11.34,-5.66)\pscircle*(15,-2){.3}
\rput(19.5,2.5){\small{$X_i$}}\rput(11.5,-1){\small{$X_j$}}\rput(17,-5){\small{$X_k$}}
\rput(22.2,-2.1){\small{$X_{ik}$}}\rput(15.1,4.8){\small{$X_{ij}$}}
\rput(10.8,-6.1){\small{$X_{jk}$}}\rput(16.3,-1.2){\small{$X_{ijk}$}}
\end{pspicture}
\caption{A 3-fold SC variety.}
\label{3conf_fig}
\end{figure}

\noindent
NC varieties often emerge as nice limits of smooth algebraic varieties.
A \sf{semistable degeneration} is a one-parameter family 
$\pi\!:\cZ\!\lra\!\De$, where $\De$ is a  disk around the origin in $\C$
and $\cZ$ is a smooth variety, such that the central fiber $\cZ_0\!\equiv\!\pi^{-1}(0)$ is an NC variety and the fibers over $\De^*\!\equiv\!\De\!-\!\{0\}$ are smooth. 
Semistable degenerations play a central role in algebraic geometry and mirror symmetry. 
They appear in compactification of moduli spaces, Hodge theory, Gromov-Witten theory, etc.

\begin{eg}[{\cite[Section~6.2]{ACGS}}]\label{QuadraticInP3_ex}
Let $P$ be a homogenous cubic polynomial in $x_0,\ldots,x_3$ and
$$\cZ'=\!\big\{\big(t,[x_0,x_1,x_2,x_3]\big)\!\in\!\C\!\times\!\P^3\!:
x_1x_2x_3=tP(x_0,x_1,x_2,x_3)\big\}\subset\C\!\times\!\P^3\,.$$
Let $\pi'\!:\cZ'\!\lra\!\C$ be the projection map to the first factor. 
For a generic $P$ and $t\!\neq\!0$, 
$\pi'^{-1}(t)$ is a smooth cubic hypersurface (divisor) in~$\P^3$. 
For $t\!=\!0$, $\pi'^{-1}(0)$ is the SC variety 
\begin{gather*}
X_{\eset}'=\{0\}\!\times\!\big(X_1'\!\cup\!X_2'\!\cup\!X_3'\big)\subset \{0\}\!\times\!\P^3
\qquad\hbox{with}\\
X_i'\!\equiv\!(x_i=0)\approx\P^2~~\forall~i\!\in\!\{1,2,3\},  \quad 
X_{ij}'\!\equiv\!X_i'\!\cap\!X_j'\approx\P^1~~\forall~i,j\!\in\!\{1,2,3\},~i\!\neq\!j.
\end{gather*}
However, the total space $\cZ'$ of $\pi'$ is not smooth at the 9 points of 
$$\cZ'^{\sing}\equiv
\{0\}\!\times\!\big(X_{\prt}'\!\cap\!(P\!=\!0)\big)\subset X_{\eset},
\qquad\hbox{where}\quad
X_{\prt}'=X_{12}'\!\cup\!X_{13}'\!\cup\!X_{23}'\subset\P^3\,.$$
A small resolution $\cZ$ of $\cZ'$ can be obtained by blowing up each singular point 
on $X_{ij}'$ in either $X_i'$ or~$X_j'$. 
The map~$\pi'$ then induces a projection $\pi\!:\cZ\!\lra\!\De$ and defines 
a semistable degeneration.
Every fiber of~$\pi$ over~$\C^*$ is a smooth cubic surface. 
The central fiber $\pi^{-1}(0)$ is the SC variety 
$X_{\eset}\!\equiv\!X_1\!\cup\!X_2\!\cup\!X_3$ with 3~smooth components, 
each a blowup of~$\P^2$ at some number of points.
If each singular point on $X_{ij}'$ is blown up in $X_i'$ with $i\!<\!j$,
then $\cZ$ is obtained from~$\cZ'$ through two global blowups of $\C\!\times\!\P^3$
and is thus projective.
\end{eg}

\noindent
As a first step to answer Gromov's inquiry, 
we introduce topological notions of NC symplectic divisor and variety
in \cite{FMZDiv,FMZNC}. 

\begin{dfn}[{\cite[Definition~2.1]{FMZDiv}}]\label{SCD_dfn}
An \sf{SC symplectic divisor} in a  symplectic manifold $(X,\om)$ 
is a finite transverse union $V\!\equiv\!\bigcup_{i \in S} V_i$ 
of smooth symplectic divisors $\{V_i\}_{i\in S}$ such that for every $I\!\subset\!S$
the submanifold
$$V_I\!\equiv\!\bigcap_{i\in I} V_i\subset X$$
is symplectic and its symplectic and intersection orientations are the same.
\end{dfn}

\noindent
Such a collection $\{V_i\}_{i\in S}$ of symplectic submanifolds of~$(X,\om)$
is sometimes called \sf{positively intersecting}.
For example, if the real dimension of~$X$ is~4,  
$V$ is an SC symplectic divisor if and only if 
the $V_i$'s intersect transversely and 
every point of a pairwise intersection $V_i\!\cap\!V_j$ is positive.
By \cite[Example~2.7]{FMZDiv}, the compatibility-of-orientation
condition along all strata cannot be replaced with a condition on a smaller subset of such strata.

\begin{dfn}[{\cite[Definition~2.5]{FMZDiv}}]\label{SCV_dfn}
An \textsf{SC symplectic variety} is a pair $(X_\eset,(\om_i)_{i \in S})$,
where
$$X_{\eset}=\bigg(\bigsqcup_{i\in S}X_i\bigg)\!\Big/\!\!\!\sim,
\quad  X_{ij}\approx X_{ji}~~\forall~i,j\!\in\!S,~i\!\neq\!j,$$
for a finite collection $(X_i,\om_i)_{i\in S}$ of symplectic manifolds,
some SC symplectic divisor $\{X_{ij}\}_{j \in S-\{i\}}$ in $X_i$ for each $i\!\in\!S$, 
and symplectic identifications $X_{ij}\!\approx\!X_{ji}$ 
for all $i,j\!\in\!S$ distinct.
\end{dfn}

\noindent
An \sf{NC symplectic divisor} is a subset of a symplectic manifold $(X,\om)$ 
locally equal to an SC symplectic divisor.
An \sf{NC symplectic variety} is a space locally equal to an SC symplectic variety. 
In other words, it is a topological space together with ``charts" 
mapping homeomorphically to SC symplectic varieties and with structure-preserving overlap maps.
There are global descriptions of NC symplectic divisors and varieties 
in terms of transverse immersions that respect certain permutation symmetries;
see \cite{FMZNC} for details.

\section{Regularizations} 
\label{RegDfn_sec}

\noindent
In order to show that Definitions~\ref{SCD_dfn} and~\ref{SCV_dfn} are appropriate analogues 
of the corresponding notions in algebraic geometry and 
are suitable for applications in symplectic topology, 
we have to show that the SC symplectic divisors and varieties of 
Definitions~\ref{SCD_dfn} and~\ref{SCV_dfn} admit (possibly after deformations)
a contractible space of compatible almost complex structures.
In the case of a smooth divisor, one needs a regularization~\eref{Reg1_e}
to construct ``nice" almost complex structures.
We need a similar thing for NC symplectic divisors and varieties.
In other words, we need to construct ``nice" neighborhoods of 
an NC symplectic divisor~$V$ in the ambient manifold~$X$ 
and of the singular locus~$X_{\prt}$ of an SC symplectic variety~$X_{\eset}$ 
(after possibly deforming these objects).
We describe below what regularizations for SC/NC divisors/varieties are.
The precise definitions are contained in \cite{FMZDiv,FMZNC}.\\

\noindent
Let $V$ be an SC symplectic divisor  in a symplectic manifold $(X,\om)$ 
as in Definition~\ref{SCD_dfn}.
For $I'\!\subset\!I\!\subset\!S$, define
$$\cN_{I;I'}=\bigoplus_{i\in I-I'}\!\!\cN_XV_i|_{V_I} \subset \cN_XV_I\,.$$
We denote by 
$$\pi_I\!:\cN_X V_I\lra V_I, \quad
\Pi_I\!: TX|_{V_I}\lra \cN_X V_I,\quad
\pi_{I;I'} : \cN_X V_I\!=\!\cN_{I;I'}\!\oplus\!\cN_{I;I-I'} \lra \cN_{I;I'}$$
the natural projection maps.

\begin{dfn}\label{RegV_dfn}
A \sf{regularization} for $V$ in $X$ is a collection of smooth embeddings
$$\Psi_I\!: \cN'_X V_I \lra X, \quad I \subset S,$$
from open neighborhoods $\cN'_XV_I\!\subset\!\cN_X V_I$ of $V_I$
so that for all $I' \!\subset\!I\!\subset\!S$
$$\Psi_I^{-1}(V_{I'})\subset\cN_{I;I'}, \quad
\Dom(\Psi_I)=\fD\Psi_{I;I'}^{-1}\big(\Dom(\Psi_{I'})\big), 
\quad  \Psi_I\!=\!\Psi_{I'}\!\circ\!\fD\Psi_{I;I'}|_{\Dom(\Psi_I)},$$
where 
\begin{gather*}
\fD\Psi_{I;I'}\!: \pi_{I;I'}^{-1}\big(\Dom(\Psi_I)\big) \lra \cN_XV_{I'},\\
\fD\Psi_{I;I'}\!\big(v_{I'},v_{I-I'}\big) \equiv 
\Pi_{I'}\bigg( \frac{\nd}{\nd t}\Psi_I(v_{I'},tv_{I-I'})\Big|_{t=0}\bigg)
\quad \forall (v_{I'},v_{I-I'}) \in \pi_{I;I'}^{-1}(\Dom(\Psi_I)),
\end{gather*}
is the derivative of $\Psi_I$ in the direction ``normal" to~$V_{I'}$.
\end{dfn}

\noindent
The $I\!=\!I'$ case in Definition~\ref{RegV_dfn} ensures that
the derivative of $\Psi_I$ along $V_I$ is the identity map.
If $|S|\!=\!1$, i.e.~$V$  is a smooth divisor, 
a regularization as in Definition~\ref{RegV_dfn} is 
a single map~$\Psi$ as in~\eref{Reg1_e}.

\begin{dfn}[{\cite[Definition~2.9]{FMZDiv}}]\label{omRegV_dfn}
An \sf{$\om$-regularization} for $V$ in $X$ consists of a choice of Hermitian structure
$(\fI_{I;i},\rho_{I;i},\na^{(I;i)})$ on $\cN_X V_i|_{V_I}$
for all $i\!\in\!I\!\subset\!S$ together with
a regularization for~$V$ in~$X$ as in Definition~\ref{RegV_dfn} so that
$$\Psi_I^*\om\!=\!\pi^*(\om|_{V_I})+\frac{1}{2} \sum_{i\in I} \nd(\rho_{I;i} \al_{\na^{(I;i)}})
\quad \forall~I\subset S,$$
and $\fD\Psi_{I;I'}$ is an isomorphism between the split Hermitian vector bundles
$\pi_{I;I'}|_{\Dom(\Psi_I) \cap \cN_{I;I'}}$ and
$\cN_X V_I|_{\Im(\Psi_I)}$.
\end{dfn}

\noindent
Since the intersection of the singular locus of an SC symplectic variety 
with each irreducible component~$X_i$ is
an SC symplectic divisor inside~$X_i$, one gets the following straightforward definition of a regularization for these objects.

\begin{dfn}[{\cite[Definition~2.15]{FMZDiv}}]\label{SCVreg_dfn}
Let $(X_\eset,(\om_i)_{i \in S})$ be an SC symplectic variety 
as in Definition~\ref{SCV_dfn}.
An \sf{$(\omega_i)_{i \in S}$-regularization} is an $\om_i$-regularization $\cR_i$
for the SC symplectic divisor $\{X_{ij}\}_{j \in S-\{i\}}$ in $X_i$
for each $i\! \in\!S$ so that the restrictions of $\cR_i$ and $\cR_j$
to $X_{ij}$ give identical $\om_i|_{X_{ij}}$-regularizations
of the SC symplectic divisor $\{X_{ijk}\}_{k\in S-\{i,j\}}$ in $X_{ij}$
for all distinct $i,j\!\in\!S$.
\end{dfn}

\noindent
Since an NC symplectic divisor $V\!\subset\!X$ is equal to 
an SC symplectic divisor $V_p\!\subset\!X$ near each point $p\!\in\!V$,
we can define an \sf{$\om$-regularization} for~$V$ to be 
an $\om$-regularization for $V_p$ for each $p\!\in\!V$ so that any two such 
$\om$-regularizations associated to $p,q\!\in\!V$ agree on $V_p\!\cap\!V_q$.
Regularizations for NC symplectic varieties are defined similarly;
see \cite{FMZNC} for details.\\

\noindent 
Unlike a smooth symplectic divisor, 
an NC symplectic divisor $V\!\subset\!X$ need not admit an $\om$-regularization. 
If an SC symplectic divisor~$V$ as in Definition~\ref{SCD_dfn} admits an $\om$-regularization,
then its smooth components~$V_i$ are in fact  $\om$-orthogonal.
On the other hand, many applications (such as symplectic constructions and
Gromov-Witten theory) care only about the deformation equivalence classes 
of the symplectic structure. 
Therefore, the alternative philosophy proposed in~\cite{FMZDiv} is~to\\

\begin{minipage}{6in}\label{DefPhil_minip}
{\it study NC symplectic divisors/varieties up to deformation equivalence and 
show that each deformation equivalence class has a subspace of 
sufficiently ``nice" representatives}.\\
\end{minipage}

\noindent
More concretely, for a transverse union $V\!\equiv\!\bigcup_{i\in S} V_i$ 
of closed real codimension~2 submanifolds of a manifold $X$, 
let $\Symp^+(X,V)$ be the space of all symplectic forms~$\om$ on $X$ such that $V$ 
is an SC symplectic divisor in~$(X,\om)$.
We also define a space of auxiliary data $\Aux(X,V)$ to be the space of pairs 
$(\om,\cR)$, where $\om\!\in\!\Symp^+(X,V)$ and $\cR$ 
is an $\om$-regularization of~$V$ in~$X$. 
Let 
$$\Pi\!:\Symp^+(X,V)\lra H^2(M;\R)$$ 
be the map sending $\om$ to its de Rham equivalence class~$[\om]$.
The following is a weaker version of the main result of \cite{FMZDiv} 
for SC symplectic divisors.

\begin{thm}[{\cite[Theorem 2.13]{FMZDiv}}]\label{SymptoAux_Thx}
Let $V\!\equiv\!\bigcup_{i\in S}V_i$ be a transverse union of closed real 
codimension~2 submanifolds of a manifold~$X$.
Then the projection~maps 
$$\pi\!:\Aux(X,V)\lra \Symp^+(X,V),\qquad
\pi\big|_{\Pi^{-1}(\al)}\!: \{\Pi\!\circ\!\pi\}^{-1}(\al) \lra \Pi^{-1}(\al), 
~~ \al\!\in\!H^2_{\tn{dR}}(M),$$
are weak homotopy equivalences.
\end{thm}

\noindent
There is a direct analogue of Theorem~\ref{SymptoAux_Thx} for SC symplectic varieties;
see \cite[Theorem 2.17]{FMZDiv}.
There are also similar results for NC symplectic divisors and varieties; see~\cite{FMZNC}.\\

\noindent
For many applications, the most important consequences of Theorem~\ref{SymptoAux_Thx} 
are the following.
First, for each $\om\!\in\!\Symp^+(X,V)$
there exists a path $(\om_t)_{t\in [0,1]}$ of cohomologous symplectic forms 
in $\Symp^+(X,V)$ such that $\om_0\!=\!\om$ and $V$ admits an $\om_1$-regularization in~$X$.
By the Moser Isotopy Theorem \cite[Section 3.2]{MS1},
this is equivalent to saying that every SC symplectic divisor in $(X,\omega)$
is isotopic through SC symplectic divisors inside $(X,\omega)$ to one which admits an $\omega$-regularization.
Second, for two pairs 
$$(\om_0,\cR_0),(\om_1,\cR_1) \in\Aux(X,V)$$ 
and a path $(\om_t)_{t\in [0,1]}$ of symplectic forms $\Symp^+(X,V)$ 
connecting~$\om_0$ and~$\om_1$, there exists 
a deformation $(\om_{t,\tau})_{t,\tau\in[0,1]}$ of this path 
in $\Symp^+(X,V)$ fixing the end points,~i.e.
$$\big(\om_{t,0}\big)_{t\in[0,1]}=\big(\om_t\big)_{t\in [0,1]}, \qquad 
\om_{0,\tau}=\om_0,~\om_{1,\tau}=\om_1~~\forall\,\tau\!\in\![0,1],$$
such that the ending path $(\om_{t,1})_{t\in [0,1]}$ can be lifted to a~path
$$\big(\om_{t,1},\wt\cR_t\big)_{t\in[0,1]}\in\Aux(X,V)
\quad\tn{with}~~\wt\cR_0\cong\cR_0,~~\wt\cR_1\cong\cR_1.$$
If $(\om_t)_{t \in [0,1]}$ is a path of cohomologous forms,
then its deformation $(\om_{t,\tau})_{t,\tau\in[0,1]}$ can be chosen to consist 
of cohomologous forms as~well.\\

\noindent
By Theorem~\ref{SymptoAux_Thx}, every SC symplectic divisor is isotopic 
to one admitting a ``nice" compatible almost complex structure. 
We define the space of almost K\"ahler data $\AK(X,V)$ to be the space 
of triples $(\om,\cR,J)$, where $\om\!\in\!\Symp^+(X,V)$ is a symplectic structure,
$\cR$ is an $\om$-regularization of $V$ as in Definition~\ref{SCVreg_dfn},
and $J\!\in\!\cJ_{\om}(X)$ is such that for each $I\!\subset\!S$:
\bIt

\item $V_I$ is $J$-holomorphic,

\item $\Psi_I^*J$ restricted to each fiber $F$ of $\pi_I|_{\Dom(\Psi_I)}$ 
is equal to $\bigoplus_{i\in I}\fI_{I;i}|_F$,  and

\item $\pi_I \circ \Psi_I^{-1}$ is a $(J|_{V_I},J|_{\Im(\Psi_I)})$-holomorphic map.

\eIt
By Theorem~\ref{SymptoAux_Thx} and a straightforward induction argument, 
the projection maps
\begin{gather*}
\pi'\!:\AK(X,V)\lra \Aux(X,V), \\
\pi'|_{\{\Pi\circ\pi\}^{-1}(\al)}\!: 
\{\Pi\!\circ\!\pi\!\circ\!\pi'\}^{-1}(\al)\lra\{\Pi\!\circ\!\pi\}^{-1}(\al), 
~~\al\!\in\!H^2_{\tn{dR}}(M),
\end{gather*}
are also weak homotopy equivalences. 
Therefore, the above conclusions concerning lifts from $\Symp^+(X,V)$
hold with $\Aux(X,V)$ replaced by $\AK(X,V)$ as well.\\

\noindent
If  $(\om,\cR,J)\!\in\!\AK(X,V)$, $J$ is very regular around $V$.
In particular, it respects the $\C^*$-action on the components $\cN_XV_i|_{V_I}$
of the normal bundle~$\cN_XV_I$ of~$V_I$ in~$X$ and
the image of its Nijenhuis tensor on $TX|_{V_I}$ lies in~$TV_I$
(i.e.~it vanishes in the normal direction to~$V_I$).
Regularizations for NC divisors and varieties similarly ensure the existence 
of almost complex structures with analogous properties on symplectic manifolds
containing NC divisors and on NC varieties themselves.
The above properties of~$J$ are very desirable for applications involving 
$J$-holomorphic curve techniques; 
works such as \cite{LR,IPrel,Jun,BP,ACGS} make use of these properties
in contexts involving various specializations of NC symplectic divisors
and varieties introduced in~\cite{FMZDiv,FMZNC}.\\

\noindent
In algebraic geometry, one can associate a holomorphic line bundle to any Cartier divisor.
It is straightforward to extend this to the symplectic topology category in
the case of a smooth symplectic divisor~$V$ in a symplectic manifold~$(X,\om)$.
Fix an identification~$\Psi$ as in~\eref{Reg1_e} and
an $\om|_{\cN_XV}$-compatible complex structure~$\fI$ on~$\cN_XV$ so that 
$\cN_XV$ becomes a complex line bundle over~$V$. 
Then,
\begin{gather*}
\cO_X(V)\equiv \big(\Psi^{-1\,*}\pi_{\cN_XV}^*\cN_XV|_{\Psi(\cN'_XV)} \sqcup 
(X\!-\!V)\!\times\!\C\big)\big/\!\!\sim\,\lra X,\\
\notag
\Psi^{-1\,*}\pi_{\cN_XV}^*\cN_XV|_{\Psi(\cN'_XV)} 
  \ni\big(\Psi(v),v,cv\big)\sim\big(\Psi(v),c\big)\in (X\!-\!V)\!\times\!\C,
\end{gather*}
is a complex line bundle over $X$ with 
\BE{c1cOXV_e}c_1(\cO_X(V))=\tn{PD}_X([V]_X)\in H^2(X;\Z),\EE
where $[V]_X$ is the homology class in~$X$ represented by~$V$.
The space of pairs~$(\Psi,\fI)$ involved in explicitly constructing this line bundle
is contractible. Therefore, the deformation equivalence class of $\cO_X(V)$  depends
only on the deformation equivalence class of~$\om$ in 
$$\Symp^+(X,V)=\Symp(X,V),$$ 
i.e.~it is an invariant of $\pi_0(\Symp^+(X,V))$.\\

\noindent
Regularizations for NC symplectic divisors and varieties extend the construction 
of the previous paragraph to these spaces.
In particular, an NC symplectic divisor~$V$ in a symplectic manifold $(X,\om)$
determines a complex line bundle $\cO_X(V)$ over~$X$ satisfying~\eref{c1cOXV_e}. 
It also determines a complex vector bundle $TX(-\!\log V)$ of rank equal
to half the real dimension of~$X$ satisfying
\BE{c1logXV_e}c\big(TX(-\!\log V)\big)=
c(TX,\om)\big/\big(1\!+\!\tn{PD}_X([V^{(1)}]_X)\!+\!\tn{PD}_X([V^{(2)}]_X)
\!+\!\ldots\big) \in H^2(X;\Q),\EE
where $V^{(r)}\!\subset\!V$ is the $r$-fold locus (locally intersection of
at least $r$ branches of~$V$); see \cite{FMZNC}.
This vector bundle extends the notion of \sf{logarithmic tangent bundle},
which plays a central role in the Gross-Siebert program \cite{GS0,GS}, 
to symplectic topology.
The deformation equivalence classes of both bundles 
depend only on the deformation equivalence class of~$\om$ in $\Symp^+(X,V)$.
If $V$ is an SC symplectic divisor as in Definition~\ref{SCD_dfn}, then
$$\cO_X(V)=\bigotimes_{i\in S}\cO_X(V_i)\lra X$$
and the equality in~\eref{c1logXV_e} holds in $H^2(X;\Z)$.
A symplectic NC variety $(X_{\eset},\om_{\eset})$ determines a complex line bundle
$\cO_{X_{\prt}}(X_{\eset})$ over the singular locus~$X_{\prt}$ of~$X_{\eset}$,
which we call the \sf{normal bundle} of~$X_{\prt}$.
The deformation equivalence class of $\cO_{X_{\prt}}(X_{\eset})$ similarly depends
only on the deformation equivalence class of~$\om_{\eset}$ in 
the space $\Symp^+(X_{\eset})$ of all NC-symplectic-variety structures on~$X_{\eset}$.
If $(X_{\eset},\om_{\eset})$ is an SC symplectic variety as in Definition~\ref{SCV_dfn}, then 
a regularization for $(X_{\eset},\om_{\eset})$ determines line bundles
$$\cO_{X_i^c}(X_i)\lra X_i^c\!\equiv\!\bigcup_{j\in S-\{i\}}\!\!\!\!\!\!X_j\subset X_{\eset},
\qquad i\in S,$$
obtained by canonically identifying the line bundles 
$$\cO_{X_j}(X_{ij})\big|_{X_{jk}}=\cO_{X_{jk}}(X_{ijk})
=\cO_{X_k}(X_{ik})\big|_{X_{jk}}$$
over $X_{jk}$; see \cite[Section~2.1]{FMZSum}.
In this case,
$$\cO_{X_{\prt}}(X_{\eset})=\bigotimes_{i\in S}\cO_{X_i^c}(X_i)\big|_{X_{\prt}}
\lra X_{\prt}\equiv\!\!
\bigcup_{\begin{subarray}{c}i,j\in S\\ i\neq j\end{subarray}}\!\!X_{ij}\,.$$

\section{Smoothings of NC symplectic varieties}
\label{SympSum_sec}

\noindent
Having introduced analogues of NC divisors/varieties into symplectic topology,
we next present an analogue of the algebro-geometric notion of
semistable degeneration. 

\begin{dfn}[{\cite[Definition 2.6]{FMZSum}},\cite{FMZNC}]\label{OneSmoothing_Df}
If $(\cZ,\om_\cZ)$ is a symplectic manifold and  $\De\!\subset\!\C$ is a disk around the origin, 
a smooth surjective map $\pi\!:\cZ\!\lra\!\De$ is 
a \sf{nearly regular symplectic fibration} if
\bIt

\item $X_{\eset}\!\equiv\!\cZ_0\!\equiv\!\pi^{-1}(0)$ is an NC symplectic divisor in $(\cZ,\om_\cZ)$, 

\item $\pi$ is a submersion outside of the singular locus $X_{\prt}$ of $X_{\eset}$,

\item for every $\la\!\in\!\De\!-\!\{0\}$, 
$\cZ_{\la}\!\equiv\!\pi^{-1}(\la)$ is a symplectic submanifold of $(\cZ,\om_\cZ)$.

\eIt
\end{dfn}

\vspace{.1in}

\noindent
The restriction of $\om_{\cZ}$ to $X_{\eset}$ above determines an NC symplectic variety
$(X_{\eset},\om_{\eset})$.
We call the nearly regular symplectic fibration of Definition~\ref{OneSmoothing_Df}
a \sf{one-parameter family of smoothings} of $(X_{\eset},\om_{\eset})$.
From the complex geometry point of view, it replaces
the nodal singularity $z_1\!\ldots\!z_N\!=\!0$ in~$\C^n$,  
i.e.~a union of $N$ coordinate hyperplanes in the central fiber~$\pi^{-1}(0)$,
by a smoothing $z_1\!\ldots\!z_N\!=\!\la$ with $\la\!\in\!\C^*$ in a generic fiber. 
A regularization for $(X_{\eset},\om_{\eset})$ determines a complex vector bundle 
$\log_{\cZ}\!TX_{\eset}$ over~$X_{\eset}$ of rank equal
to half the real dimension of~$X_{\eset}$ satisfying
$$c\big(\log_{\cZ}\!TX_{\eset}\big)=c\big(T\cZ(-\!\log X_{\eset})\big)\big|_{X_{\eset}}\,.$$
This bundle is the analogue of \sf{logarithmic tangent bundle}
of a smoothable NC variety in algebraic geometry.
Its deformation equivalence class again depends only
on the deformation equivalence class of~$\om_{\cZ}$ in $\Symp^+(\cZ,X_{\eset})$.\\

\noindent
In light of the deformation equivalence philosophy stated on page~\pageref{DefPhil_minip},
we say that an NC symplectic variety $(X_{\eset},\om_{\eset})$ is \sf{smoothable} 
if some NC symplectic variety $(X_{\eset},\om'_{\eset})$ deformation equivalent 
to $(X_{\eset},\om_{\eset})$  admits a one-parameter family of smoothings.
In this section, we answer the following question:\\

\begin{minipage}{6in}
{\it Which SC symplectic varieties are smoothable?}\\
\end{minipage}

\noindent
The \sf{$d$-semistability condition} of \cite[Definition (1.13)]{F} is well-known 
to be an obstruction to the smoothability of an NC variety in a one-parameter family 
with a smooth total space in the algebraic geometry category. 
As shown in~\cite{PP}, the $d$-semistability condition is not 
the only obstruction in the algebraic category, even in 
the simplest non-trivial case discussed below.\\

\noindent
Let $(X_1,\om_1)$ and $(X_2,\om_2)$ be smooth symplectic manifolds 
with identical copies of a smooth symplectic divisor $X_{12}\!\subset\!X_1,X_2$. 
Then 
\BE{2foldSC_e}
\big(X_\eset\!=\!X_1\!\cup_{X_{12}}\!X_2, \om_\eset\!=\!(\om_1,\om_2)\big)\EE
is an SC symplectic variety.
If  $(X_1,\om_1)$, $(X_2,\om_2)$, and $X_{12}$ are smooth projective varieties,
the $d$-semistability condition of~\cite{F} in this case is the existence
of an isomorphism
\BE{Obundle_e}
\cO_{X_{\prt}}(X_{\eset})\!\equiv\!
\cN_{X_1}X_{12}\otimes_{\C} \cN_{X_2}X_{12} \approx \cO_{X_{12}}\lra X_{12}\EE
in the category of {\it holomorphic} line bundles.
By the now classical \sf{symplectic sum construction}, suggested in~\cite[p343]{Gr} 
and carried out in~\cite{Gf,MW}, the smoothability of
the 2-fold SC symplectic variety~\eref{2foldSC_e} in the sense of Definition~\ref{OneSmoothing_Df}
is {\it equivalent} to
the existence of an isomorphism~\eref{Obundle_e} 
in the category of {\it complex} line bundles.\\

\noindent
The topological type of smooth fibers $(X_{\#},\om_{\#})$ of the one-parameter family 
of smoothings produced by the construction of~\cite{Gf} 
depends only on  the homotopy class of isomorphisms~\eref{Obundle_e}.
With such a choice fixed, this construction involves choosing 
an $\om_1|_{\cN_{X_1}X_{12}}$-compatible almost complex structure on~$\cN_{X_1}X_{12}$,
an $\om_2|_{\cN_{X_2}X_{12}}$-compatible almost complex structure on~$\cN_{X_2}X_{12}$,
and a representative for the above homotopy class.
Because of these choices, the resulting symplectic manifold $(X_{\#},\om_{\#})$
is determined by $(X_1,\om_1)$, $(X_2,\om_2)$, and the choice of the homotopy class
only up to symplectic deformation equivalence.
Since the symplectic deformations of the SC symplectic variety (\ref{2foldSC_e})
do not affect  the deformation equivalence class of $(X_{\#},\om_{\#})$,
it would have been sufficient to carry out the symplectic sum construction
of~\cite{Gf} only on a path-connected set of representatives for
each deformation equivalence class of the SC symplectic variety (\ref{2foldSC_e}).\\

\noindent
The above change in perspective turns out to be very useful for smoothing out 
arbitrary NC symplectic varieties in \cite{FMZSum,FMZNC} and
thus answering another question of  \cite[p343]{Gr}.
By the next theorem, the direct analogue of the $d$-semistability condition of~\cite{F}
is the only obstruction for the smoothability of an arbitrary NC symplectic variety.
Regularizations are the essential auxiliary data in the proof of this result.

\begin{thm}[{\cite[Theorem~2.7]{FMZSum}},\cite{FMZNC}]\label{SympSum_thm}
An NC symplectic variety $(X_\eset,\om_{\eset})$ is smoothable if and only~if 
the associated line bundle $\cO_{X_{\prt}}(X_{\eset})$ is trivializable. 
Furthermore, the germ at the zero fiber of the deformation equivalence class of 
the nearly regular symplectic fibration $(\cZ,\om_{\cZ},\pi)$ provided 
by the proof of this statement is determined by a homotopy class of trivializations 
of~$\cO_{X_{\prt}}(X_{\eset})$.
If in addition $X_{\prt}$ is compact,  
the deformation equivalence class of a smooth fiber~$(\cZ_{\la},\om_{\la})$
is also determined by a homotopy class of these trivializations.
\end{thm}

\noindent
If $(X_{\eset},\om_{\eset})$ is an SC symplectic variety as in Definition~\ref{SCV_dfn}
 and $X_{\eset}$ is compact,
we call (the deformation equivalence class of) a generic fiber 
of the resulting one-parameter family 
the \sf{multifold} or \sf{$|S|$-fold symplectic sum} of $(X_i,\om_i)_{i\in S}$.\\

\noindent
If $(X_{\eset},\om_{\eset})$ is an SC symplectic variety as in Definition~\ref{SCV_dfn},
\BE{ResToI_e}\cO_{X_{\prt}}(X_{\eset})\big|_{X_{ij}} = 
\cN_{X_i}X_{ij}\otimes\cN_{X_j}X_{ij}\otimes
\bigotimes_{\begin{subarray}{c}k\in S\\ k\neq i,j\end{subarray}} 
\!\!\cO_{X_{ij}}(X_{ijk}) \qquad\forall~i,j\!\in\!S,\,i\!\neq\!j.\EE
By \cite[Example~2.10]{FMZSum}, the triviality of $\cO_{X_{\prt}}(X_{\eset})$ 
is in general stronger than the triviality of the restrictions~\eref{ResToI_e};
the latter is known as the \sf{triple point condition} in the algebraic geometric literature.
If the natural homomorphism
\BE{H2restr_e}
H^2(X_{\prt};\Z) \lra 
\bigoplus_{\begin{subarray}{c}i,j \in S\\ i\neq j\end{subarray}}\!\!H^2(X_{ij};\Z)\EE
is injective or at most one of the homomorphisms
\BE{H2restr_e2}
H^1(X_{ij};\Z)\lra 
\bigoplus_{\begin{subarray}{c}k \in S\\ k\neq i,j\end{subarray}}\!\!H^1(X_{ijk};\Z),
\qquad i,j\!\in\!S,\,i\!\neq\!j,\EE
is not surjective, then $\cO_{X_\prt}(X_\eset)$ is trivial if and only if 
all restrictions in~\eref{ResToI_e} are trivial.\\

\noindent
In Example~\ref{QuadraticInP3_ex}, the line bundle (\ref{ResToI_e}) corresponding 
to $X_{ij}'\!\approx\!\P^1$ is equal to~$\cO_{\P^1}(3)$.
Therefore, $X_{\eset}'$ is not smoothable.
After blowing up the 3 singular points on each~$X_{ij}'$, we~get  
$$\cN_{X_{i}}X_{ij}\otimes \cN_{X_{j}}X_{ij}\otimes \cO_{X_{ij}}(X_{123})\!\cong\!\cO_{\P^1}
\qquad\forall~i,j\!=\!1,2,3,~i\!\neq\!j.$$
In this case, the homomorphism~\eref{H2restr_e} is injective 
and all homomorphisms~\eref{H2restr_e2} are surjective.
By Theorem~\ref{SympSum_thm},
the NC symplectic variety~$X_{\eset}$ of Example~\ref{QuadraticInP3_ex} is thus smoothable.
In this case, there is only one homotopy class of trivializations of~$\cO_{X_{\prt}}(X_{\eset})$.
The smoothing provided by the proof of Theorem~\ref{SympSum_thm} is conjecturally equivalent
to the one of Example~\ref{QuadraticInP3_ex}, provided the latter is projective.\\

\noindent
The proof of the SC case of Theorem~\ref{SympSum_thm} in~\cite{FMZSum}
explicitly constructs $(\cZ,\om_{\cZ},\pi)$ by gluing together local charts~$\cZ_I$
with $I\!\subset\!S$ non-empty;
these charts are equipped with symplectic forms and local smoothings.
We deform $(\om_i)_{i \in S}$ in $\Symp^+(X_{\eset})$
so that $X_{\eset}$ admits an $(\om_i)_{i \in S}$-regularization~$\cR$
and choose a trivialization~$\Phi$ of $\cO_{X_{\prt}}(X_{\eset})$ so that
it is compatible with~$\cR$ in a suitable sense.
If $|I|\!\ge\!2$, $\cZ_I$ is a neighborhood of a large open subset~$X_{I}^{\circ}$
of~$X_I$ in 
$$\cN X_I|_{X_I^{\circ}}\equiv \bigoplus_{i\in I} \cN_{X_{I-i}} X_I|_{X_I^\circ};$$
the closure of $X_I^{\circ}$ is ``slightly" disjoint from all $X_J$ with $J\!\supsetneq\!I$.
If $I\!=\!\{i\}$ with $i\!\in\!S$,
$\cZ_I$ is a neighborhood of $X_i^{\circ}\!\times\!\{0\}$ in $X_i^\circ\!\times\!\C$.  
The regularization~$\cR$ and the trivialization~$\Phi$ are used to glue 
the charts~$\cZ_{I'}$ and~$\cZ_I$ for $I'\!\subset\!I$ 
with $|I'|\!\ge\!2$ and $|I'|\!=\!1$, respectively.\\

\noindent
If $|I|\!\ge\!2$, the restriction of~$\pi$ to~$\cZ_I$ is a positive multiple 
of the composition
$$\cZ_I \lhra  \bigoplus_{i\in I}\cN_{X_{I-i}}X_I\big|_{X_I^\circ}
\lra  \bigotimes_{i\in I}\cN_{X_{I-i}}X_I\big|_{X_I^\circ}
=\cO_{X_\partial}(X_\eset)|_{X_I^\circ}
\stackrel{\Phi}{\lra} X_I^\circ\!\times\!\C\lra \C.$$
If $|I|\!=\!1$, the restriction of~$\pi$ to~$\cZ_I$ is a positive multiple 
of the projection to the second component.
The symplectic form~$\om_{\cZ}$ on~$\cZ$ is built by interpolating between 
the symplectic forms~\eref{OmtoHatOm_e} determined by the product Hermitian structures
$(\fI_{I;i},\rho_{I;i},\na^{(I;i)})_{i\in I}$ on $\cZ_I$ with $|I|\!\ge\!2$
and certain product symplectic forms on  $\cZ_I$ with $|I|\!=\!1$.\\

\noindent
The proof of the SC case of Theorem~\ref{SympSum_thm} outlined above is extended
to the general case in~\cite{FMZNC}.
By \cite[Proposition~5.1]{FMZSum} and its extension to the NC case in~\cite{FMZNC}, 
every one-parameter family of smoothings of an NC symplectic variety $(X_{\eset},\om_{\eset})$
determines a homotopy class of trivializations of $\cO_{X_\partial}(X_\eset)$. 
By \cite[Proposition~5.5]{FMZSum} and its extension to the NC case in~\cite{FMZNC},
the homotopy class determined by the family provided by the proof of Theorem~\ref{SympSum_thm}
is the input homotopy class.
We believe that the equivalence classes of smoothings of a compact NC symplectic variety
$(X_{\eset},\om_{\eset})$ correspond to 
the homotopy classes of trivializations of $\cO_{X_{\prt}}(X_{\eset})$.
This is equivalent to the following.

\begin{cnj}\label{smoothing_cnj}
Let $(\cZ,\om_{\cZ},\pi)$ be a one-parameter family of smoothings
of a compact NC symplectic variety $(X_{\eset},\om_{\eset})$.
Then $(\cZ,\om_{\cZ},\pi)$ is deformation equivalent to a smoothing
of $(X_{\eset},\om_{\eset})$ provided by the proof of Theorem~\ref{SympSum_thm}.
\end{cnj}

\begin{rmk}
The surgery construction of~\cite{SymingtonThesis,Symington3} 
on 4-dimensional symplectic manifolds along  pairwise positively intersecting immersed surfaces, 
also called $N$-fold symplectic sum construction, agrees with ours 
(which is consistent with algebraic geometry and \cite[p343]{Gr}) only for $N\!=\!3$.
In particular, the setting of \cite[Theorem~2.7]{Symington3} is essentially
the $\dim_{\R}X\!=\!4$ case of the setting of \cite[Theorem~2.7]{FMZSum}.
The output of \cite[Theorem~2.7]{Symington3} is then symplectically deformation equivalent
to the smooth fibers of the one-parameter family provided by \cite[Theorem~2.7]{FMZSum}.
The perspectives taken in~\cite{Symington3} and in~\cite{FMZSum} are fundamentally different
as~well. 
The viewpoint in~\cite{Symington3} is that of surgery on 4-dimensional manifolds;
the viewpoint in~\cite{FMZSum} is that of smoothing a variety in a one-dimensional family 
with  a smooth total space.
The configurations in~\cite{Symington3} with $N\!\ge\!4$ correspond to varieties, such~as
\BE{Gfconf_e4}
\big\{(x,y,z,w)\!\in\!\C^4\colon\!\,xy\!=\!0,\,zw\!=\!0\big\},\EE
that do not even admit such smoothings.
The total space of the natural one-parameter smoothing of~(\ref{Gfconf_e4}),
i.e.~with~$0$ replaced by~$\la\!\in\!\C$, is singular at the origin.
\end{rmk}

\section{SC degenerations of symplectic manifolds} 
\label{SympCut_sec}

\noindent
We next discuss the potential for reversing the construction of Theorem~\ref{SympSum_thm}.\\

\begin{minipage}{6in}
\textit{Can one degenerate a symplectic manifold $(X,\om)$ into 
some SC symplectic variety $(X_i,\om_i)_{i\in S}$ in a one-parameter family?}\\
\end{minipage}

\noindent
The $|S|\!=\!2$ case of this question is the now classical 
\sf{symplectic cut construction} of~\cite{L}.
Given a free Hamiltonian $S^1$-action generated by Hamiltonian~$h$ on an~open subset~$W$ of~$X$ 
so that $\wt{V}\!\equiv\!h^{-1}(0)$ is a separating hypersurface, 
this construction decomposes~$(X,\om)$ into two symplectic manifolds,
$(X_-,\om_-)$ and~$(X_+,\om_+)$.
It cuts~$X$ into closed subsets~$U^{\le}$ and~$U^{\ge}$ 
along~$\wt{V}$ and collapses their boundary~$\wt{V}$
to a smooth symplectic divisor $V\!\equiv\!\wt{V}/S^1$ inside~$(X_-,\om_-)$ and~$(X_+,\om_+)$.
The associated ``wedge"  
$$X_{\eset}\equiv X_-\!\cup_{V}\!X_+$$
is a 2-fold SC~symplectic variety as in~\eref{2foldSC_e}.
If we assume instead that $\wt{V}$ is non-separating, 
the result would be an NC symplectic variety.\\

\noindent
In \cite{FZCut}, we generalize and enhance the symplectic cut construction of~\cite{L}
to produce a nearly regular symplectic fibration with regular fibers 
deformation equivalent to $(X,\om)$.
The central singular fiber of this fibration is what we call 
a \sf{multifold (or $N$-fold) symplectic cut} of $(X,\om)$. 
The input of this construction is a multifold cutting configuration}~$\sC$ 
defined below.\\

\noindent
For a finite non-empty set $S$, let 
$$(S^1)^S_{\bu}=\bigg\{\!(\ne^{\fI\th_i})_{i\in S}\!\in\!(S^1)^S\!\!:
\prod_{i\in S}\!\ne^{\fI\th_i}\!=\!1\!\bigg\}.$$
For $I\!\subset\!S$, we identify $(S^1)^I$ with the subgroup
$$\big\{(\ne^{\fI\th_i})_{i\in S}\!\in\!(S^1)^S\!: \ne^{\fI\th_i}\!=\!1~\forall\,
i\!\in\!S\!-\!I\big\}$$
of $(S^1)^S$ in the natural way and let
$$(S^1)^I_{\bu}\!\equiv\!(S^1)^S_{\bu}\!\cap\!(S^1)^I\,.$$
Denote by $\ft_{I;\bu}\!\subset\!\ft_{S;\bu}$ the Lie algebra of $(S^1)^I_{\bu}$
and by~$\ft_{I;\bu}^*$ its dual.
For $i,j\!\in\!I\!\subset\!S$, the homomorphism
$$\ft_{I;\bu}^*\!=\!\R^I\big/\big\{a^I\!\in\!\R^I\!:\,a\!\in\!\R\big\}\lra\R,
\qquad \eta\!\equiv\!\big[(a_k)_{k\in I}\big]\lra \eta_{ij}\!\equiv\!a_j\!-\!a_i\,,$$
is well-defined.
We write $(\eta)_i\!<\!(\eta)_j$ (resp.~$(\eta)_i\!\le\!(\eta)_j$, $(\eta)_i\!=\!(\eta)_j$) 
if $0\!<\!\eta_{ij}$  (resp.~$0\!\le\!\eta_{ij}$, $0\!=\!\eta_{ij}$).

\begin{dfn}\label{SympCut_dfn0}
A \sf{multifold Hamiltonian configuration} for a symplectic manifold~$(X,\om)$ is a~tuple 
\BE{SympCutDfn_e}
\sC\equiv\big(U_I,\phi_I,\mu_I\!:X\!\lra\!\ft_{I;\bu}^*\big)_{\eset\neq I\subset S},\EE
where $S$ is a finite non-empty set,
$(U_I)_{\eset\neq I\subset S}$ is an open cover of~$X$, and 
$\phi_I$   is a Hamiltonian $(S^1)^I_{\bu}$-action on~$U_I$ with moment map~$\mu_I$, 
such~that 
\begin{enumerate}[label=$(\alph*)$,leftmargin=*]

\item\label{IJinter_it} $U_I\!\cap\!U_{I'}\!=\!\eset$ unless $I\!\subset\!I'$ or $I'\!\subset\!I$;

\item\label{phiIUJ_it} $\mu_I(x)|_{\ft_{I';\bu}}\!=\!\mu_{I'}(x)$
for all $x\!\in\!U_I\!\cap\!U_{I'}$ and $I'\!\subset\!I\!\subset\!S$;

\item\label{UIJpos_it} $(\mu_I(x))_i\!<\!(\mu_I(x))_j$ 
for all $x\!\in\!U_I\!\cap\!U_{I'}$, $i\!\in\!I'\!\subset\!I\!\subset\!S$, 
and $j\!\in\!I\!-\!I'$.
\end{enumerate}
\end{dfn}

\begin{dfn}\label{SympCut_dfn1}
A \sf{multifold cutting configuration for~$(X,\om)$} is 
a  multifold Hamiltonian configuration  as in~\eref{SympCutDfn_e}
such~that the restriction of the $(S^1)^I_{\bu}$-action~$\phi_I$ to $(S^1)^{I'}_{\bu}$ 
is free on the preimage of $0\!\in\!\ft_{I';\bu}^*$ under the moment~map
$$\mu_{I';I}\!: \big\{x\!\in\!U_I\!\!:
\big(\mu_I(x)\!\big)_{\!i}\!<\!\big(\mu_I(x)\!\big)_{\!j}
~\forall\,i\!\in\!I',\,j\!\in\!I\!-\!I'\big\}\lra \ft_{I';\bu}^*,~~
\mu_{I';I}(x)\!=\!\mu_I(x)\big|_{\ft_{I';\bu}}\,,$$
for all \hbox{$I'\!\subset\!I\!\subset\!S$} with $I'\!\neq\!\eset$.
\end{dfn}

\noindent
We use a multifold cutting configuration~$\sC$ in~\cite{FZCut} 
to decompose~$(X,\om)$ into $|S|$~symplectic manifolds $(X_i,\om_i)$ at~once.
We first cut~$X$ into the closed subspaces
$$U_i^{\leq}\equiv \bigcup_{i\in I\subset S}\!\!\!\!
\big\{x\!\in\!U_I\!:\,(\mu_I(x))_i\!\le\!(\mu_I(x))_j~\forall\,j\!\in\!I\big\},
\qquad i\!\in\!S.$$ 
The subspace~$U_i^{\leq}$ has boundary and corners
$$U_I^{\leq} \equiv  \bigcup_{I\subset J\subset S}\!\!\!\!
\big\{x\!\in\!U_J\!:\,(\mu_J(x))_i\!\le\!(\mu_J(x))_j~\forall\,i\!\in\!I,~j\!\in\!J\big\},
\qquad i\!\in\!I\!\subset\!S,~|I|\!\ge\!2.$$ 
We collapse each~$U_I^{\leq}$ by the $(S^1)^I_{\bu}$-action~$\phi_I$  to
obtain symplectic manifolds $(X_i,\om_i)$ with $i\!\in\!S$ and
 symplectic submanifolds $X_I\!\equiv\!U_I^{\leq}/(S^1)^I_{\bu}$ 
of real codimension $2(|I|\!-\!1)$ in $(X_i,\om_i)$ with $i\!\in\!I\!\subset\!S$.
For each $i\!\in\!S$, $\{X_{ij}\}_{i\in S-i}$ is an SC symplectic divisor in $(X_i,\om_i)$.
The entire collection $\{X_I\}$ determines an SC~symplectic variety
$(X_{\eset},\om_{\eset})$, which we call
a \sf{multifold (or $|S|$-fold) symplectic cut of $(X,\om)$}.\\

\noindent
We also construct a symplectic manifold~$(\cZ,\om_{\cZ})$ containing $X_{\eset}$ 
as an SC symplectic divisor and a smooth map $\pi\!:\cZ\!\lra\!\C$ 
so that $X_{\eset}\!=\!\pi^{-1}(0)$ and $\om_{\cZ}|_{X_{\eset}}\!=\!\om_{\eset}$.
The restriction of~$\pi$ to a neighborhood~$\cZ'$ of~$X_{\eset}$ in~$\cZ$
is a nearly regular symplectic fibration and thus determines
a one-parameter family of smoothings of the SC symplectic variety $(X_{\eset},\om_{\eset})$.
If $X$ is compact, then a generic fiber of $\pi|_{\cZ'}$ 
is deformation equivalent to~$(X,\om)$.
In such a case, we call $(\cZ',\om_{\cZ}|_{\cZ'},\pi|_{\cZ'})$
a \sf{multifold} (or \sf{$|S|$-fold}) \sf{SC symplectic degeneration} of~$(X,\om)$.
This is a symplectic topology analogue of the algebro-geometric notion of 
\sf{semistable degeneration} (i.e.~smooth one-parameter family of degenerations 
of a smooth algebraic variety to an NC algebraic variety).\\

\noindent
The symplectic manifold $(\cZ,\om_\cZ)$ is obtained by gluing together charts
$(\cZ_{I}^\circ,\vp_I)$ with $I\!\subset\!S$ non-empty.
Each $(\cZ_{I}^\circ,\vp_I)$ is the Symplectic Reduction \cite[Theorem~23.1]{daSilva} 
with respect to the moment~map
$$\wt\mu_I\!:U_I\!\times\!\C^I\lra\ft_{I;\bu}^*, \qquad
\wt\mu_I(x,z)=\mu_I(x)-\mu_{\C^I;\bu}(z),$$
where $\mu_{\C^I;\bu}$ is the moment map for the restriction of the standard 
$(S^1)^I$-action on $\C^I$ to $(S^1)^I_{\bu}$
such that $0\!\in\!\C^I$ is its critical point.
The restriction of~$\pi$ to~$\cZ_{I}^\circ$ is a positive multiple of
the~map
$$\cZ_I^{\circ}\lra\C, \qquad 
\big[x,(z_i)_{i\in I}\big] \lra\prod_{i\in I}z_i\,.$$
The intersection of $\cZ_I^\circ$ with $X_i$ corresponds to 
the region $z_i\!=\!0$ inside $\cZ_I^\circ$.\\

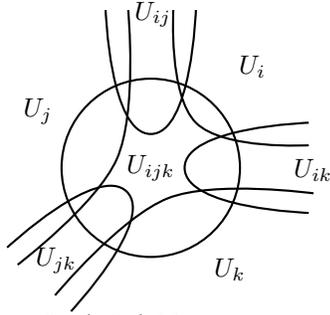
\begin{figure}
\begin{pspicture}(-3,-2)(11,1)
\psset{unit=.3cm}
\pscurve(22,-4)(16.5,-2)(22,0)\rput(22.2,-2.1){\small{$U_{ik}$}}
\psrotate(15,-2){90}{\pscurve(22,-4)(16.5,-2)(22,0)}\rput(15.1,4.8){\small{$U_{ij}$}}
\psrotate(15,-2){225}{\pscurve(22,-4)(16.5,-2)(22,0)}\rput(10.8,-6.1){\small{$U_{jk}$}}
\pscircle(15,-2){4}\rput(15,-2){\small{$U_{ijk}$}}
\pscurve(16,5)(17,0)(22,-1)\rput(19.5,2.5){\small{$U_i$}}
\pscurve(9.1,-6.3)(13.5,-1.6)(14,5)\rput(10,.5){\small{$U_j$}}
\psrotate(15,-2){135}{\pscurve(9.1,-6.3)(13.5,-1.6)(14,5)}\rput(18.5,-6.5){\small{$U_k$}}
\end{pspicture}
\caption{A 3-fold cutting configuration.}
\label{3cut_fig}
\end{figure}

\noindent
The $S\!=\!\{1,2\}$ case of Definition~\ref{SympCut_dfn1} corresponds to 
the symplectic cut construction of~\cite{L}
by identifying $(S^1)^S_{\bu}$ with~$S^1$ via the projection to
the first component of~$(S^1)^S$ and taking 
$$W=U_S, \quad X_-=U_1^{\le}, \quad X_+=U_2^{\le}, \quad 
\wt{V}=U_{12}^{\le}\,.$$
Figures~\ref{3cut_fig} and~\ref{3conf_fig} show 
a 3-fold cutting configuration and the associated 3-fold symplectic cut,
respectively.\\

\noindent
The SC case of the symplectic sum/smoothing construction of Theorem~\ref{SympSum_thm} 
and the symplectic cut/degeneration construction of~\cite{FZCut} are intuitively 
mutual inverses.
However, we are unaware of any work where even the 2-fold case of
this statement (relating the constructions of~\cite{Gf} and~\cite{L})
is made precise.
The purpose of~\cite{FZRelation} is to establish this statement as formulated below.\\

\noindent
Fix $n,N\!\in\!\Z^+$.
Let $\tn{SCV}(n,N)$ be the space of tuples $(X_{\eset},\om_{\eset},\hb)$  
consisting of a compact $N$-fold SC symplectic variety $(X_{\eset},\om_{\eset})$ 
of real dimension~$2n$ 
and a homotopy class~$\hb$ of trivializations of an associated line bundle 
$\cO_{X_\prt}(X_{\eset})$.
Let $\tn{SCC}(n,N)$ be the space of tuples $(X,\om,\sC)$ consisting 
of a compact symplectic manifold $(X,\om)$  
of real dimension~$2n$ and an $N$-fold cutting configuration~$\sC$ for~$(X,\om)$.
By \cite{FZRelation}, a generic fiber $(X,\om)$ of 
a nearly regular symplectic fibration arising from the proof of
the SC case of Theorem~\ref{SympSum_thm} admits a natural $N$-fold cutting configuration~$\sC$.
By \cite[Proposition~5.1]{FMZSum}, the semistable degeneration 
$(\cZ',\om_{\cZ}|_{\cZ'},\pi|_{\cZ'})$ of~$(X,\om)$ arising from a cutting configuration~$\sC$
determines 
a homotopy class~$\hb$ of trivialization of a line bundle $\cO_{X_\prt}(X_{\eset})$
associated with the central fiber.
Thus, there are natural~maps
\BE{cScC_e}\cS_{n,N} \colon\! \tn{SCV}(n,N) \lra \tn{SCC}(n,N), 
\qquad \cC_{n,N} \colon\! \tn{SCC}(n,N) \lra \tn{SCV}(n,N),\EE
which we call \sf{smoothing} and \sf{cutting} maps respectively.\\

\noindent
The aim of \cite{FZRelation} is to show that the two maps in~\eref{cScC_e}
are weak homotopy inverses.
It is fairly straightforward to show that $\cC_{n,N}\!\circ\!\cS_{n,N}$
is homotopy equivalent to the identity.
Conjecture~\ref{smoothing_cnj} essentially implies that  $\cS_{n,N}\!\circ\!\cC_{n,N}$
is weakly homotopy equivalent to the identity, but is more than what is needed.
The multifold symplectic cut construction of~\cite{FZCut}, 
the maps~\eref{cScC_e}, and their being weak homotopy inverses
should generalize to the arbitrary NC case as~well.

\section{Directions for further research} 
\label{FurRes_sec}

\noindent
\sf{Log smooth degenerations} to \sf{log smooth} algebraic varieties
play important roles in such areas of modern algebraic geometry
as Hodge theory and \sf{log Gromov-Witten theory} \cite{GS,AC}.
The almost K\"ahler analogue of the log smooth category provided by 
the \sf{exploded manifold} category of \cite{BP-EX} underpins 
a similar study of GW-invariants in~\cite{BP}. 
The works~\cite{GS,AC,BP} extend GW-invariants to
algebraic varieties with so-called \sf{fine saturated} log structures
and show that these invariants do not change under 
deformations that are smooth in the category of such varieties.
An effective decomposition formula splitting GW-invariants of 
a fine saturated log algebraic variety into GW-invariants of 
the irreducible components of the underlying variety
would  generalize 
the renowned formula of~\cite{Jun} and is key to the Gross-Siebert program~\cite{GS0},
but has turned out to be difficult to work~out.
As GW-invariants of smooth algebraic varieties are fundamentally symplectic topology invariants, 
it is natural to expect the same of fine saturated log algebraic varieties.
This should in turn provide a more robust setting for effectively generalizing 
the decomposition formula of~\cite{Jun}.
We believe that the deformation equivalence philosophy stated on page~\pageref{DefPhil_minip}
and the methods used to implement it in the case of NC singularities can be used
to define a category of \sf{fine saturated log symplectic varieties},
which in turn should provide a suitable setting for extending GW-invariants
from symplectic manifolds and studying their properties under semistable degeneration.
NC~singularities are the most basic type of fine saturated log structures.
Our work \cite{FMZDiv,FMZSum,FMZNC,FZCut,FZRelation} introduces symplectic topology analogues
of these structures, including the associated \sf{log tangent bundles}, and  lays
the foundation for defining the long-awaited GW-invariants relative to NC
symplectic~divisors.\\

\noindent
The  ($2$-fold) symplectic sum construction of~\cite{Gf} has been used to
build vast classes of non-K\"ahler symplectic manifolds.
For example, it is shown in~\cite{Gf} that every finitely presented group
can be realized as the fundamental group of a compact symplectic manifold of real dimension~4. 
There are still many open questions about the geography of symplectic manifolds 
and the classification of symplectic manifolds with certain topological properties.
The smoothing and degeneration constructions of \cite{FMZSum,FMZNC,FZCut,FZRelation}
may shed light on some of these questions.
For example, it is well-known that a \sf{rationally connected} compact K\"ahler manifold
(i.e.~a K\"ahler manifold with a rational curve through every pair of points) is simply connected;
see \cite[Theorem~3.5]{Campana}.
As noted by J.~Starr, the fundamental group of a compact almost K\"ahler manifold $(X,\om,J)$
with a rational $J$-holomorphic curve of a fixed homology class
through every pair of points is finite.
The existence of such curves is implied by  the existence of a nonzero GW-invariant
of~$(X,\om)$ with two point insertions, 
but the converse is not known to be true even in the projective category.

\begin{ques}\label{RCsymp_ques}
Is every compact almost K\"ahler manifold with a rational 
$J$-holomorphic curve of a fixed homology class through every pair of points simply connected?
Is every compact symplectic manifold $(X,\om)$ with  a nonzero GW-invariant
with two point insertions simply connected?
\end{ques}

\noindent
The multifold sum/smoothing construction of~\cite{FMZSum} 
and the multifold cut/degenerative construction of~\cite{FZCut} 
may be useful in answering these questions negatively and positively,
respectively.\\

\noindent
{\it Simons Center for Geometry and Physics, Stony Brook University, Stony Brook, NY 11794\\
mtehrani@scgp.stonybrook.edu}\\

\noindent
{\it Department of Mathematics, Stony Brook University, Stony Brook, NY 11794\\
markmclean@math.stonybrook.edu, azinger@math.stonybrook.edu}


\begin{thebibliography}{99}

\bibitem{AC} D.~Abramovich and Q.~Chen, 
{\it Stable logarithmic maps to Deligne-Faltings pairs~II},
Asian J.~Math.~18 (2014), no.~3, 465--488

\bibitem{ACGS} D.~Abramovich, Q.~Chen, M.~Gross, and B.~Siebert, 
{\it Decomposition of degenerate Gromov-Witten invariants},
preprint available from the authors

\bibitem{Campana} F.~Campana, 
{\it On twistor spaces of the class $\sC$},
J.~Diff.~Geom.~33 (1991), no.~2, 541-–549

\bibitem{daSilva} A.~Canas da Silva, 
{\it Lectures on Symplectic Geometry},
Lecture Notes in Mathematics~1764, Springer-Verlag, 2001 (revised~2006)

\bibitem{Donaldson} S.~Donaldson, 
{\it  Symplectic submanifolds and almost-complex geometry}, 
J.~Differential Geom.~44 (1996), no.~4, 666--705

\bibitem{FMZDiv} M.~Farajzadeh-Tehrani, M.~McLean, and A.~Zinger,
{\it Normal crossings singularities for symplectic topology}, math/1410.0609

\bibitem{FMZSum} M.~Farajzadeh-Tehrani, M.~McLean, and A.~Zinger,
{\it The smoothability of normal crossings symplectic varieties}, math/1410.2573

\bibitem{FMZNC} M.~Farajzadeh-Tehrani, M.~McLean, and A.~Zinger,
{\it Towards a theory of singular symplectic varieties}, in preparation

\bibitem{FZCut} M.~Farajzadeh-Tehrani and A.~Zinger,
{\it Normal crossings degenerations of symplectic manifolds},
math/1603.0766

\bibitem{FZRelation}
M.~Farajzadeh Tehrani and A. Zinger, 
{\it On the multifold symplectic sum and cut constructions}, work in progress

\bibitem{F} R.~Friedman,
{\it Global smoothings of varieties with normal crossings},
Ann.~Math.~118 (1983), no.~1, 75--114

\bibitem{Gf} R.~Gompf, {\it A new construction of symplectic manifolds}, 
Ann.~of Math.~142 (1995), no.~3, 527--595

\bibitem{Gr} M.~Gromov, {\it Partial Differential Relations}, Springer-Verlag, 1986

\bibitem{GS0} M.~Gross and B.~Siebert,
{\it  Affine manifolds, log structures, and mirror symmetry},
Turkish J.~Math.~27 (2003), no.~1, 33-60

\bibitem{GS}
M.~Gross and B.~Siebert, 
{\it Logarithmic Gromov-Witten invariants}, JAMS~26 (2013), no.~2, 451--510

\bibitem{IPrel} E.~Ionel and T.~Parker, 
{\it Relative Gromov-Witten invariants}, Ann.~of Math.~157 (2003), no.~1, 45--96. 

\bibitem{L} E.~Lerman, {\it Symplectic cuts}, 
Math.~Res.~Lett.~2 (1995), no.~3, 247--258

\bibitem{LR}  A.-M.~Li and Y.~Ruan, 
{\it Symplectic surgery and Gromov-Witten invariants of Calabi-Yau 3- folds}, 
Invent.~Math.~145 (2001), no.~1, 151--218

\bibitem{Jun} J.~Li, {\it A degeneration formula for GW-invariants},
J.~Diff.~Geom.~60 (2002), no.~1, 199--293

\bibitem{LT}  J.~Li and G.~Tian, 
{\it Virtual moduli cycles and Gromov-Witten invariants of general symplectic manifolds}, 
Topics in Symplectic \hbox{$4$-Manifolds},
47-83, First Int.~Press Lect.~Ser., I, Internat.~Press, 1998

\bibitem{MW} J.~McCarthy and J.~Wolfson, 
{\it Symplectic normal connect sum}, Topology 33 (1994), no.~4, 729--764

\bibitem{MS1} D.~McDuff and D.~Salamon, 
{\it Introduction to Symplectic Topology}, 2nd Ed., Oxford University Press, 1998

\bibitem{MS2} D.~McDuff and D.~Salamon, 
{\it J-Holomorphic Curves and Symplectic Topology}, 
Colloquium Publications 52, AMS, 2012

\bibitem{MAff} M.~McLean, 
{\it The growth rate of symplectic homology and affine varieties}, 
Geom.~Funct. Anal.~22 (2012), no.~2, 369--442

\bibitem{McLean} M.~McLean, 
{\it Reeb orbits and the minimal discrepancy of an isolated singularity},
Invent.~Math.~204 (2016), no.~2, 505--594

\bibitem{Morrison} D.~Morrison,
{\it Compactifications of moduli spaces inspired by mirror symmetry}, math/9304007

\bibitem{NN}  A.~Newlander and L.~Nirenberg,
{\it Complex analytic coordinates in almost complex manifolds}, 
Ann.~Math.~65 (1957), no.~3, 391--404

\bibitem{BP-EX}  B.~Parker, {\it Exploded fibrations}, 
Proceedings of Gokova Geometry-Topology Conference 2006 (2007), 52--90

\bibitem{BP} B.~Parker,
{\it Gromov-Witten invariants of exploded manifolds},
math/1102.0158

\bibitem{BP3} B.~Parker,
{\it Gluing formula for Gromov-Witten invariants in a triple product},
math/1511.00779

\bibitem{PP} U.~Persson and H.~Pinkham, 
{\it Some examples of non-smoothable varieties with normal crossings}, 
Duke Math.~J.~50 (1983), no.~2, 477--486

\bibitem{Sheridan} N.~Sheridan,
{\it Homological mirror symmetry for Calabi-Yau hypersurfaces in projective space},
Invent.~Math.~199 (2015), no.~1, 1--186

\bibitem{SymingtonThesis} M.~Symington, 
{\it New constructions of symplectic four-manifolds}, Ph.D.~Thesis, Stanford University,~1996

\bibitem{Symington3} M.~Symington,  {\it A new symplectic surgery: the 3-fold sum},
 Topology and its Applications 88 (1998), no.~1, 27--53

\bibitem{Tian}
G. Tian,  {\it The quantum cohomology and its associativity,} 
Current Developments in Mathematics (1995), 361--401, Inter.~Press

\bibitem{TsengYau} L.-S.~Tseng and S.-T.~Yau.
{\it Non-K\"{a}hler Calabi-Yau manifolds}, 
String-Math 2011, Proc.~Sympos.~Pure Math.~85 (2011), 241--254

\end{thebibliography}
\end{document}